\documentclass[12pt]{article}
\usepackage{amsmath,amsfonts,amssymb,amsthm}
\begin{document}

\newcommand{\1}{{{\bf 1}}}
\newcommand{\id}{{\rm id}}
\newcommand{\Hom}{{\rm Hom}\,}
\newcommand{\End}{{\rm End}\,}
\newcommand{\Res}{{\rm Res}\,}
\newcommand{\Image}{{\rm Im}\,}
\newcommand{\Ind}{{\rm Ind}\,}
\newcommand{\Aut}{{\rm Aut}\,}
\newcommand{\Ker}{{\rm Ker}\,}
\newcommand{\gr}{{\rm gr}}
\newcommand{\Der}{{\rm Der}\,}
\newcommand{\Z}{\mathbb{Z}}
\newcommand{\Q}{\mathbb{Q}}
\newcommand{\C}{\mathbb{C}}
\newcommand{\N}{\mathbb{N}}
\newcommand{\g}{\mathfrak{g}}
\newcommand{\gl}{\mathfrak{gl}}
\newcommand{\h}{\mathfrak{h}}
\newcommand{\wt}{{\rm wt}\,}
\newcommand{\A}{\mathcal{A}}
\newcommand{\D}{\mathcal{D}}
\newcommand{\Lie}{\mathcal{L}}
\newcommand{\E}{\mathcal{E}}

\def \b{\beta}
\def \<{\langle}
\def \>{\rangle}
\def \be{\begin{equation}\label}
\def \ee{\end{equation}}
\def \bex{\begin{exa}\label}
\def \eex{\end{exa}}
\def \bl{\begin{lem}\label}
\def \el{\end{lem}}
\def \bt{\begin{thm}\label}
\def \et{\end{thm}}
\def \bp{\begin{prop}\label}
\def \ep{\end{prop}}
\def \br{\begin{rem}\label}
\def \er{\end{rem}}
\def \bc{\begin{coro}\label}
\def \ec{\end{coro}}
\def \bd{\begin{de}\label}
\def \ed{\end{de}}

\newtheorem{thm}{Theorem}[section]
\newtheorem{prop}[thm]{Proposition}
\newtheorem{coro}[thm]{Corollary}
\newtheorem{conj}[thm]{Conjecture}
\newtheorem{exa}[thm]{Example}
\newtheorem{lem}[thm]{Lemma}
\newtheorem{rem}[thm]{Remark}
\newtheorem{de}[thm]{Definition}
\newtheorem{hy}[thm]{Hypothesis}
\makeatletter \@addtoreset{equation}{section}
\def\theequation{\thesection.\arabic{equation}}
\makeatother \makeatletter

\begin{Large}
\begin{center}
\textbf{Rationality of the vertex algebra $V_L^+$ when $L$ is a nondegenerate even lattice of arbitrary rank}
\end{center}
\end{Large}
\begin{center}{
Gaywalee Yamskulna\footnote{E-mail: gyamsku@ilstu.edu}\\
Department of Mathematical Sciences, Illinois State University, Normal, IL 61790\\
and\\
Institute of Science, Walailak University, Nakon Si Thammarat,
Thailand}
\end{center}
\begin{abstract}

In this paper we prove  that the vertex algebra $V_L^+$ is rational  if $L$ is a negative definite even lattice of  finite rank , or if  $L$ is a non-degenerate even lattice of a finite rank that is neither positive definite nor negative definite. 
In particular, for such even lattices $L$, we show that the Zhu algebras of the vertex algebras $V_L^+$ are semisimple. This extends the result of Abe from \cite{a} which establishes the rationality of $V_L^+$   when $L$ is a positive definite even lattice of  finite rank.

\end{abstract}
\section{Introduction}
The vertex algebras $V_L^+$ are one of the most  important classes of vertex algebras along with those vertex algebras  associated with lattices, affine Lie algebras and Virasoro algebras. They were originally introduced in the Frenkel-Lepowsky-Meurman construction of the moonshine module vertex algebra (see \cite{flm}). The study of the representation theory of $V_L^+$ began with the case when $L$ is a positive definite even lattice. In this case, the classification of all irreducible $V_L^+$-modules, and the study of the complete reducibility property of $V_L^+$-modules was done by Abe, Dong, Jiang, and Nagatomo (see \cite{a1, a, ad, dj, dn1}). When $L$ is a rank one negative definite even lattice, the classification of irreducible $V_L^+$-modules was completed by Jordan in \cite{j}. Later, in \cite{y}, the author classified all irreducible $V_L^+$-modules for the case when $L$ is a negative definite even lattice of arbitrary rank, and when $L$ is a non-degenerate even lattice that is neither positive definite nor negative definite. 
 
In this paper, we continue our study of the representation theory of $V_L^+$ when $L$ is a negative definite even lattice of  finite rank and when $L$ is a non-degenerate even lattice that is neither positive definite nor negative definite. We prove here  that $V_L^+$ is rational in these cases. The main idea of the proof is to show that the Zhu algebras $A(V_L^+)$ are semisimple. Note that a vertex algebra $V$ is called rational if any $V$-module is completely reducible.

This paper is organized as follows. In Section 2, we review the necessary background material. In particular, we recall the construction of irreducible $V_L^+$-modules and recall certain facts about the Zhu algebras $A(V_L^+)$ when $L$ is a negative definite even lattice of a finite rank and when $L$ is a non-degenerate even lattice that is neither positive definite nor negative definite. The proof of the rationality of $V_L^+$ in these cases is contained in
section 3.
\section{Preliminaries}

First, we discuss relationships between vertex algebras and Zhu algebras. Next, we recall the constructions of irreducible $V_L^+$- modules, and review the Zhu algebras $A(V_L^+)$ when $L$ is a negative definite even lattice of a finite rank and when $L$ is a non-degenerate even lattice that is neither positive definite nor negative definite.
\subsection{Relationships between Vertex Algebras and Zhu algebras}
\begin{de}\cite{lli} A {\em vertex algebra} $V$ is a vector space equipped with a linear map
$Y(\cdot, z):V\rightarrow(\End V)[[z,z^{-1}]], v\mapsto Y(v,z)=\sum_{n\in\Z}v_nz^{-n-1}$ and a distinguished vector ${\bf 1}\in V$ which satisfies the following properties:
\begin{enumerate}
\item $u_nv=0$ for $n>>0$.
\item $Y({\bf 1},z)=id_V$.
\item $Y(v,z){\bf 1}\in V[[z]]$ and $\lim_{z\rightarrow 0}Y(v,z){\bf 1}=v$.
\item (the Jacobi identity) \begin{eqnarray*}
& &z_0^{-1}\delta\left(\frac{z_1-z_2}{z_0}\right)Y(u,z_1)Y(v,z_2)-z_0^{-1}\delta\left(\frac{z_2-z_1}{-z_0}\right)Y(v,z_2)Y(u,z_1)\\
&&=z_2^{-1}\delta\left(\frac{z_1-z_0}{z_2}\right)Y(Y(u,z_0)v,z_2).
\end{eqnarray*}
\end{enumerate}
We denote the vertex algebra just defined by $(V, Y, {\bf 1})$ or, briefly, by $V$.
\end{de}
\begin{de} A $\Z$-graded vertex algebra is a vertex algebra $$V=\oplus_{n\in\Z}V_n; \text{ for }v\in V_n,\ \  n=\wt v,$$ equipped with a conformal vector $\omega\in V_2$ which satisfies the following relations:
\begin{itemize}
\item $[L(m),L(n)]=(m-n)L(m+n)+\frac{1}{12}(m^3-m)\delta_{m+n,0}c_V$ for $m,n\in \Z$, where $c_V\in \C$ (the central charge) and $$Y(\omega,z)=\sum_{n\in\Z}L(n)z^{-n-2}\left(=\sum_{m\in\Z}\omega_mz^{-m-1}\right);$$
\item $L(0)v=nv=(\wt v) v$ for $n\in\Z$, and $v\in V_n$;
\item $Y(L(-1)v,z)=\frac{d}{dz}Y(v,z)$.
\end{itemize}
\end{de}

For the rest of this subsection, {\em we assume that $V$ is a $\Z$-graded vertex algebra}.
\begin{de}\cite{dlm2} A {\em weak $V$-module} $M$ is a vector space equipped with a linear map $Y_M(\cdot, z):V\rightarrow(\End  M)[[z]], v\mapsto Y_M(v,z)=\sum_{n\in\Z}v_nz^{-n-1}$ which satisfies the following properties: for $v,u\in V$, and $w\in M$
\begin{enumerate}
\item $v_nw=0$ for $n>>0$.
\item $Y_M({\bf 1},z)=id_M$.
\item (the Jacobi identity) \begin{eqnarray*}
& &z_0^{-1}\delta\left(\frac{z_1-z_2}{z_0}\right)Y_M(u,z_1)Y_M(v,z_2)-z_0^{-1}\delta\left(\frac{z_2-z_1}{-z_0}\right)Y_M(v,z_2)Y_M(u,z_1)\\
&&=z_2^{-1}\delta\left(\frac{z_1-z_0}{z_2}\right)Y_M(Y(u,z_0)v,z_2).
\end{eqnarray*}
\end{enumerate}
\end{de}

\begin{de} An {\em irreducible} weak $V$-module is a weak $V$-module that has no weak $V$-submodule except 0 and itself. Here, a weak sub-module is defined in the obvious way.
\end{de}

\begin{de}\cite{dlm2} An ({\em ordinary}) $V$-module is a weak $V$-module $M$ which carries a $\C$-grading induced by the spectrum of $L(0)$. Then $M=\oplus_{\lambda\in\C}M_{\lambda}$ where $M_{\lambda}=\{w\in M| L(0)w=\lambda w\}$, and $\dim M_{\lambda} <\infty$. Moreover, for fixed $\lambda$, $M_{n+\lambda}=0$ for all small enough integers $n$.
\end{de}
\begin{de}\cite{dlm2} An {\em admissible} $V$-module $M$ is a $\Z_{\geq 0}$-graded weak $V$-module $M=\oplus_{n\in\Z_{\geq 0}}M(n)$ such that $v_mM(n)\subset M(n+\wt v-m-1)$ for any homogeneous $v\in V$ and $m,n\in\Z$. Here, $\Z_{\geq 0}$ is the set of nonnegative integers.

An admissible $V$-submodule of $M$ is a weak $V$-submodule $N$ of $M$ such that $N=\oplus_{n\in\Z_{\geq 0}}N\cap M(n)$.
\end{de}

\begin{de} An {\em irreducible admissible} $V$-module is an admissible $V$-module that has no admissible submodule except 0 and itself.
\end{de}
\begin{de} A vertex algebra $V$ is called a {\em rational} if every admissible $V$-module is completely reducible, i.e., a direct sum of irreducible admissible $V$-modules.
.\end{de}
\begin{prop}\cite{dlm2, z}

\begin{enumerate}
\item Any ordinary $V$-module is an admissible $V$-module.
\item For any irreducible admissible $V$-module $M$, there exists a complex number $\lambda$ such that $M=\oplus_{n=0}^{\infty}M(\lambda+n)$ where $M(\lambda+n)$ is the $L(0)$-eigenspace of the eigenvalue $\lambda+n$. We call $\lambda$ the lowest weight of $M$.
\end{enumerate}
\end{prop}

Next, we will define a Zhu algebra and we will discuss the relationships between the Zhu algebra and a vertex algebra. 

For a homogeneous vector $u\in V$, $v\in V$, we define products $u*v$, and $u\circ v$ as follow:
\begin{eqnarray*}
u*v&=&\Res_{z}\left(\frac{(1+z)^{\wt u}}{z} Y(u,z)v\right)\\
u\circ v&=&\Res_{z}\left(\frac{(1+z)^{\wt u}}{z^2}Y(u,z)v\right).
\end{eqnarray*}  Then we extend these products linearly on $V$. We let $O(V)$ be the linear span of $u\circ v$ for all $u,v\in V$ and we set $A(V)=V/O(V)$. Also, for $v\in V$, we denote $v+O(V)$ by $[v]$.
\begin{thm}\label{zhu}\cite{z}\ \
\begin{enumerate}
\item $(A(V),*)$ is an associative algebra with the identity $[{\bf 1}]$. Moreover, $[\omega]$ is a central element of $A(V)$.
\item The map $M\rightarrow M(0)$ gives a bijection between the set of equivalence classes of irreducible admissible $V$-modules to the set of equivalence classes of simple $A(V)$-modules.
\end{enumerate}
\end{thm}

We denote by $P(V)$ the set of lowest weights of all irreducible admissible $V$ -modules. The following is a key proposition.
\begin{prop}\label{pv}\cite{a1} If the Zhu algebra $A(V)$ is
semisimple, and $(\lambda+\Z_+)\cap P(V)=\emptyset$ for any $\lambda\in P(V)$, then $V$ is rational. Here, $\Z_+$ is the set of positive integers.
\end{prop}

\subsection{Vertex algebras $V_L^+$, $M(1)^+$ and their Zhu algebras}
First, we will discuss the construction of the vertex algebras $M(1)^+$, $V_L^+$ and irreducible $V_L^+$-modules. Next, we will review some information about the Zhu algebras $A(M(1)^+)$ and $A(V_L^+)$.

We will follow the setting in \cite{flm}. Let $L$ be a non-degenerate even lattice of a rank $d$. We set $\hat{L}$ be the canonical central extension of $L$ by the cyclic group $\<\kappa\>$ of order 2. We let $e:L\rightarrow\hat{L}$ be a section such that $e_0=1$ and we let $\epsilon:L\times L\rightarrow\<\kappa\>$ be the corresponding 2-cocycle. We assume that $\epsilon$ is bimultiplicative. Then $\epsilon(\alpha,\beta)/\epsilon(\beta, \alpha)=\kappa^{\<\alpha,\beta\>}$, $\epsilon(\alpha,\beta)\epsilon(\alpha+\beta, \gamma)=\epsilon(\beta,\gamma)\epsilon(\alpha, \beta+\gamma),$ and $e_{\alpha}e_{\beta}=\epsilon(\alpha,\beta)e_{\alpha+\beta}$ for $\alpha,\beta, \gamma\in L$.  Let $\theta$ denote an automorphism of $\hat{L}$ defined by $\theta(e_{\alpha})=e_{-\alpha}$ and $\theta(\kappa)=\kappa$. Furthermore, we set $K=\{a^{-1}\theta(a)|a\in\hat{L}\}$. 

We define $V_L=M(1)\otimes\C[L]$ where $M(1)$ is the Heisenberg vertex operator algebra associated with $\h=\C\otimes_{\Z}L$ and $\C[L]$ is the group algebra of $L$ with basis vectors $e^{\alpha}$ for $\alpha\in L$. Note that $\C[L]$ is an $\hat{L}$-module under the action $e_{\alpha}e^{\beta}=\epsilon(\alpha,\beta)e^{\alpha+\beta}$. It was shown in \cite{b,flm, gu} that $V_L$ is a $\Z$-graded simple vertex algebra. Moreover, $M(1)$ is its $\Z$-graded vertex sub-algebra.

Next, we define a linear automorphism $\theta:V_{L}\rightarrow V_{L}$ by
$$\theta(\b_1(-n_1)\b_2(-n_2)....\b_k(-n_k)e^{\alpha})=(-1)^k\b_1(-n_1)....\b_k(-n_k)e^{-\alpha}.$$ Consequently, $\theta Y(v,z)u=Y(\theta v,z)\theta(u)$ for $u,v\in V_L$. In particular, $\theta$ is an automorphism of $V_L$ and $M(1)$.

For any stable $\theta$-subspace $U$ of $V_{L}$, we denote a $\pm 1$ eigenspace of $U$ for $\theta$ by $U^{\pm}$.
\begin{prop}\cite{dm} $M(1)^+$ and $V_L^+$ are simple $\Z$-graded vertex algebras.
\end{prop}
\begin{rem} $V_L^+$ is not an admissible module of itself when $L$ is not a positive definite even lattice.
\end{rem}


We let $\h[-1]$ and $M(1)(\theta)$ be the $\theta$-twisted Heisenberg algebra and its unique irreducible module, respectively. We set $\chi$ be a central character of $\hat{L}/K$ such that $\chi(\kappa)=-1$ and we let $T_{\chi}$ be the irreducible $\hat{L}/K$-module with central character $\chi$.  
We define $V_L^{T_{\chi}}=M(1)(\theta)\otimes T_{\chi}.$ Note that $V_L^{T_{\chi}}$ is an irreducible $\theta$-twisted $V_L$-module (cf. \cite{d, flm}).

We define an action of $\theta$ on $M(1)(\theta)$ and $V_L^{T_{\chi}}$ in the following way:
\begin{eqnarray*}
& &\theta(\b_1(-n_1)\b_2(-n_2)....\b_k(-n_k){\bf 1})=(-1)^k\b_1(-n_1)....\b_k(-n_k){\bf 1},\text{ and }\\
& &\theta(\b_1(-n_1)\b_2(-n_2)....\b_k(-n_k)t)=(-1)^k\b_1(-n_1)....\b_k(-n_k)t).
\end{eqnarray*}
for $\beta_i\in\h$, $n_i\in\frac{1}{2}+\Z_{\geq 0}$ and $t\in T_{\chi}$. We denote by $M(1)(\theta)^{\pm}$ and $V_L^{T_{\chi},{\pm}}$ the $\pm 1$-eigenspace for $\theta$ of $M(1)(\theta)$ and $V_L^{T_{\chi}}$, respectively.
\begin{thm}\label{y}\cite{y} Let $L$ be an even lattice of a finite rank. Suppose that $L$ is either a negative definite lattice or a non-degenerate lattice that is neither positive definite nor negative definite. Then the set of all irreducible admissible $V_L^+$-modules is 
\begin{eqnarray*}
\{\ \ V_L^{T_{\chi}, \pm}&|& T_{\chi} \text{ is irreducible } \hat{L}/K-\text{module with central character } \chi\\& & \text{ such that }\chi(\kappa)=-1\ \ \}\end{eqnarray*}  \end{thm} 
\begin{coro}\label{pvl+} Assume that $L$ is a rank $d$ non-degenerate even lattice. If $L$ is either a negative definite lattice or a non-degenerate lattice that is neither positive definite nor negative definite, then $P(V_L^+)=\{\frac{d}{16},\frac{d+8}{16} \}$.
\end{coro}
Next, we recall the Zhu algebras of $M(1)^+$ and $V_L^+$ when $L$ is a negative definite even lattice and when $L$ is a non-degenerate even lattice that is neither positive definite nor negative definite. 

We let $L$ be a rank $d$ even lattice with a nondegenerate symmetric $\Z$-bilinear form $\<\cdot,\cdot\>$. We set $\h=\C\otimes_{\Z}L$ and extend $\<\cdot,\cdot\>$ to a $\C$-bilinear form on $\h$. 
Let $\{ h_a|1\leq a\leq d\}$ be an orthonormal basis of $\h$, and set $\omega_a=\omega_{h_a}=\frac{1}{2}h_a(-1)^2{\bf 1}$ and $J_a=h_a(-1)^4{\bf 1}-2h_a(-3)h_a(-1){\bf 1}+\frac{3}{2}h_a(-2)^2{\bf 1}$.  Note that vectors $\omega_a$, and $J_a$ generate a vertex operator algebra $M(1)^+$ associated to the one-dimensional vector space $\C h_a$. 

Following \cite{ad}, we set $S_{ab}(m,n)=h_a(-m)h_b(-n){\bf 1}$, and define $E^u_{ab}$, $E^t_{ab}$, and $\Lambda_{ab}$ as follows:
\begin{eqnarray*}
E^u_{ab}&=&5S_{ab}(1,2)+25S_{ab}(1,3)+36S_{ab}(1,4)+16S_{ab}(1,5)\ \ (a\neq b),\\
E^u_{aa}&=&E^u_{ab}*E^u_{ba},\\
E^t_{ab}&=&-16(3S_{ab}(1,2)+14S_{ab}(1,3)+19S_{ab}(1,4)+8S_{ab}(1,5))\ \ (a\neq b),\\
E^t_{aa}&=&E^t_{ab}*E^t_{ba},\\
\Lambda_{ab}&=&45S_{ab}(1,2)+190S_{ab}(1,3)+240S_{ab}(1,4)+96S_{ab}(1,5).
\end{eqnarray*}
\begin{prop}\cite{dn} For any $1\leq i,j,k,l\leq d$ we have $[E^t_{ij}]*[E^t_{kl}]=\delta_{jk}[E^t_{il}]$.\end{prop}

Let $A^u$, and $A^t$ be the linear subspaces of $A(M(1)^+)$ spanned by $[E^u_{ab}]$ and $[E^t_{ab}]$, respectively. Here, $1\leq a,b\leq d$.
\begin{prop}\label{am1gen}\cite{dn} \ \ 
\begin{enumerate}
\item The spaces $A^u$ and $A^t$ are two sided ideals of $A(M(1)^+)$. Moreover, ideals $A^u$, $A^t$, the unit $I^u=\sum_{i=1}^d[E^u_{ii}]$ of $A^u$ and the unit $I^t=\sum_{i=1}^d[E^t_{ii}]$ of $A^t$ are independent of the choice of an orthonormal basis.
\item There are algebra isomorphisms between $A^u$ and $\End M(1)^-(0)$, and between $A^t$ and $\End M(1)(\theta)^-$, respectively. In particular, under the basis $$\{h_1(-1){\bf 1},...,h_d(-1){\bf 1}\}$$ of $M(1)^-(0)$, each $[E^u_{ab}]$ corresponds to the matrix element $E_{ab}$ whose $(a,b)$-entry is 1 and zero elsewhere. Similarly, under the basis $$\{h_1(-\frac{1}{2}){\bf 1},...,h_d(-\frac{1}{2}){\bf 1}\}$$ of $M(1)(\theta)^-(0)$, each $[E^t_{ab}]$ corresponds to the matrix element $E_{ab}$ whose $(a,b)$-entry is 1 and zero elsewhere.
\item The Zhu algebra $A(M(1)^+)$ is generated by $[\omega_a]$, $[J_a]$ for $1\leq a\leq d$, $[\Lambda_{ab}]$ for $1\leq a\neq b\leq d$ and $[E^u_{ab}]$, $[E^t_{ab}]$ for $1\leq a,b\leq d$. 
\item The quotient algebra $A(M(1)^+)/(A^t+A^u)$ is commutative. Furthermore, it is generated by the images of $[\omega_a]$, $[J_a]$ for $1\leq a\leq d$ and $[\Lambda_{ab}]$ for $1\leq a\neq b\leq d$.
\end{enumerate}\end{prop}

For any $\alpha\in L$, we set $$V_L^+[\alpha]=M(1)^+\otimes E^{\alpha}\oplus M(1)^-\otimes F^{\alpha}$$ 
and $$A(V_L^+)(\alpha)=(V_L^+[\alpha]+O(V_L^+))/O(V_L^+).$$ Here, $E^{\alpha}=e^{\alpha}+e^{-\alpha}$, and $F^{\alpha}=e^{\alpha}-e^{-\alpha}$. Note that $A(V_L^+)$ is the sum of $A(V_L^+)(\alpha)$ for all $\alpha\in L$. For a $\Z$-graded vertex sub-algebra $U$ of $V_L^+$, the identity map induces an algebra homomorphism from $A(U)$ to $A(V_L^+)$. For $u\in U$, we use $[u]$ to denote $u+O(U)$ and $u+O(V_L^+)$.

Next, we recall several results in \cite{y}.
\begin{prop}\label{info1}\cite{y} Let $L$ be a negative definite even lattice of a finite rank $d$. In $A(V_L^+)$, we have the following.
\begin{enumerate}
\item For any indices $a,b$, we have $[E^u_{ab}]=0$ and $[\Lambda_{ab}]=0$.
\item For $\alpha\in L-\{0\}$, we let $\{h_1,...,h_d\}$ be an orthonormal basis of $\h$  such that $h_1\in \C\alpha$. Then 
\begin{enumerate}
\item $[E^t_{ab}]*[E^{\alpha}]=[E^{\alpha}]*[E^t_{ab}]$ if $a\neq 1$ and $b\neq 1$.
\item $[E^t_{1b}]*[E^{\alpha}]=-\frac{1}{2\<\alpha,\alpha\>-1}[E^{\alpha}]*[E^t_{1b}]$ if $b\neq 1$.
\item $[E^t_{b1}]*[E^{\alpha}]=-(2\<\alpha,\alpha\>-1)[E^{\alpha}]*[E^t_{b1}]$ if $b\neq 1$.
\item $[E^t_{aa}]*[E^{\alpha}]=[E^{\alpha}]*[E^t_{aa}]$ for $a\in\{1,...,d\}$.
\end{enumerate}
\item Let $I^t$ be the unit of the simple algebra $A^t$. Then for any $\alpha\in L$, $I^t*[E^{\alpha}]=[E^{\alpha}]*I^t$.
\item For any $\alpha\in L$, we have 
\begin{eqnarray*}
A(V_L^+)(\alpha)&=&span_{\C}\{[u]*[E^{\alpha}]| u\in M(1)^+\}\\
&=&span_{\C}\{[E^{\alpha}]*[u]| u\in M(1)^+\}.\end{eqnarray*}
\item Let $\alpha\in L-\{0\}$. We set $\{h_1,...,h_d\}$ be an orthonormal basis of $\h$ such that $h_1\in\C\alpha$. Then $$[H_a]=[H_1]-\frac{9}{8}[E^t_{aa}]+\frac{9}{8}[E^t_{11}]\text{\ \ \ \ \ \   for }a\in\{2,...,d\} .$$
 \end{enumerate}
\end{prop}
\begin{prop}\label{info2}\cite{y} Let $L$ be a non-degenerate rank $d$ even lattice that is neither positive definite nor negative definite. In $A(V_L^+)$, we have the following.
\begin{enumerate}
\item For any indices $a,b$, we have $[E^u_{ab}]=0$ and $[\Lambda_{ab}]=0$.
\item For $\alpha\in L$ such that $\<\alpha,\alpha\>\neq 0$, we let $\{h_1,...,h_d\}$ be an orthonormal basis of $\h$  such that $h_1\in \C\alpha$. Then 
\begin{enumerate}
\item $[E^t_{ab}]*[E^{\alpha}]=[E^{\alpha}]*[E^t_{ab}]$ if $a\neq 1$ and $b\neq 1$.
\item $[E^t_{1b}]*[E^{\alpha}]=-\frac{1}{2\<\alpha,\alpha\>-1}[E^{\alpha}]*[E^t_{1b}]$ if $b\neq 1$.
\item $[E^t_{b1}]*[E^{\alpha}]=-(2\<\alpha,\alpha\>-1)[E^{\alpha}]*[E^t_{b1}]$ if $b\neq 1$.
\item $[E^t_{aa}]*[E^{\alpha}]=[E^{\alpha}]*[E^t_{aa}]$ for $a\in\{1,...,d\}$.
\end{enumerate}
\item Let $I^t$ be the unit of the simple algebra $A^t$. Then for any $\alpha\in L$ such that $\<\alpha,\alpha\>\neq 0$, we have $I^t*[E^{\alpha}]=[E^{\alpha}]*I^t$.
\item Let $\alpha\in L$ such that $\<\alpha,\alpha\>\neq 0$. We set $\{h_1,...,h_d\}$ be an orthonormal basis of $\h$ such that $h_1\in\C\alpha$. Then $$[H_a]=[H_1]-\frac{9}{8}[E^t_{aa}]+\frac{9}{8}[E^t_{11}]\text{\ \ \ \ \ \   for }a\in\{2,...,d\} .$$
\item For any $\alpha\in L$ such that $\<\alpha,\alpha\>\neq 0$, we have 
\begin{eqnarray*}
A(V_L^+)(\alpha)&=&span_{\C}\{[u]*[E^{\alpha}]| u\in M(1)^+\}\\
&=&span_{\C}\{[E^{\alpha}]*[u]| u\in M(1)^+\}.\end{eqnarray*}
\item Let $\alpha\in L-\{0\}$ such that $\<\alpha,\alpha\>=0$. Then there exists $\beta\in L$ such that $\<\beta,\beta\><0$, $\<\alpha,\beta\><0$, and $A(V_L^+)(\alpha)\subset A(V_L^+)(\alpha+2\beta)$.
\item $A(V_L^+)$ is spanned by $A(V_L^+)(0)$, and $A(V_L^+)(\alpha)$ for all $\alpha\in L$ such that $\<\alpha,\alpha\>\neq 0$.
 \end{enumerate}
\end{prop}

\section{Rationality of the vertex algebra $V_L^+$ when $L$ is a nondegenerate even lattice of arbitrary rank}

Let $L$ be a nondegenerate even lattice of a finite rank $d$. If $L$ is positive definite, then it was shown in \cite{a1, a} that $V_L^+$ is a rational vertex algebra. In this section, we will extend this result to other cases.
\subsection{Case I: $L$ is a negative definite even lattice}

In this subsection, we will prove that the vertex algebra $V_L^+$ is rational when $L$ is a negative definite even lattice. The key idea is to show that the Zhu algebra $A(V_L^+)$ is semisimple. 

For the rest of this subsection, {\em we assume that $L$ is a negative definite even lattice.}  Following \cite{dj}, we set 
$$[\tilde{B}_{\alpha}]=2^{\<\alpha,\alpha\>-1}( I^t*[E^{\alpha}]-\frac{2\<\alpha,\alpha\>}{2\<\alpha,\alpha\>-1}[E^t_{11}]*[E^{\alpha}])$$  for $\alpha\in L-\{0\},$ and $[\tilde{B}_0]=I^t$. Here, $E^t_{11}$ is defined with respect to an orthonormal basis $\{h_a|1\leq a\leq d\}$ of $\h$ such that $h_1\in \C\alpha$. Clearly, for $i\in\{1,...,d\}$, we have
\begin{enumerate}
\item $[E^t_{ij}]*[\tilde{B}_{\alpha}]=2^{\<\alpha,\alpha\>-1}[E^t_{ij}]*[E^{\alpha}]$ when $j\in\{2,...,d\}$, and 
\item $[E^t_{i1}]*[\tilde{B}_{\alpha}]=-\frac{2^{\<\alpha,\alpha\>-1}}{2\<\alpha,\alpha\>-1}[E^t_{i1}]*[E^{\alpha}]$.
\end{enumerate} 
\begin{lem}\label{eb}
For $1\leq a,b\leq d$, $\alpha\in L$, $[E^t_{ab}]*[\tilde{B}_{\alpha}]=[\tilde{B}_{\alpha}]*[E^t_{ab}]$.
\end{lem}
\begin{proof} This follows immediately from Proposition \ref{info1}.\end{proof}
Next, we let $A^t_L=Span_{\C}\{[E^t_{ab}]*[\tilde{B}_{\alpha}]\ \ |\ \  1\leq a,b\leq d, \ \ \alpha\in L\}.$ We then have the following.
\begin{thm}\label{atl} \ \ 

\begin{enumerate}
\item $A^t_L$ is a 2-sided ideal of $A(V_L^+)$.
\item $A^t_L$ is a semisimple associative algebra. In particular, $A^t_L$ is isomorphic to $A^t\otimes_{\C}\C[\hat{L}/K]/J$ where $\C[\hat{L}/K]$ is the group algebra of $\hat{L}/K$ and $J$ is the ideal of $\C[\hat{L}/K]$ generated by $\kappa K+1$.
\end{enumerate}
\end{thm}
\begin{proof} 1. follows immediately from Propositions \ref{am1gen}, \ref{info1} and Lemma \ref{eb}. For 2., we will show that $A^t_L\cong A^t\otimes \C[\hat{L}/K]/J$ where $J$ is the ideal of $\C[\hat{L}/K]$ generated by $\kappa K+1$. Clearly, $A^t_L$ is an $A^t$-module and $A^t_L=A^t\cdot A^t_L$. Moreover, $A^t_L$ is a direct sum of $M(1)(\theta)^-$. By following the proof of Lemma 3.15 of \cite{y}, we then have that $[\tilde{B}_{\alpha}]*[\tilde{B}_{\beta}]=\epsilon(\alpha,\beta)[\tilde{B}_{\alpha+\beta}]$ for $\alpha,\beta\in L$. Notice that a linear map $\psi: \hat{L}\rightarrow A^t_L$ defined by $\psi(e_{\alpha})=[\tilde{B}_{\alpha}]$ and $\psi(\kappa)=-I^t$ induced a linear map $\bar{\psi}:\hat{L}/K\rightarrow A^t_L$ since $\theta(e_{\alpha})=e_{-\alpha}$ and $\tilde{B}_{\alpha}=\tilde{B}_{-\alpha}$. Moreover $\psi$ induced an injective algebra homomorphism from $\C[\hat{L}/K]/J$ into $A^t_L$ and an algebra isomorphism from $A^t\otimes \C[\hat{L}/K]/J$ onto $A^t_L$. Consequently, $A^t_L$ is semisimple.\end{proof}
\begin{coro}\label{cirrmod} The set of all irreducible $A^t_L$-modules is \begin{eqnarray*}
\{ \ \ \h(-1/2)\otimes T_{\chi}\ \ &|& \ \ T_{\chi} \text{ are irreducibles }\hat{L}/K-\text{modules}\\
& &\text{ with central character }\chi\text{ such that }\chi(\iota(\kappa))=-1\ \ \}.\end{eqnarray*} Here, $\h(-1/2)=\{h(-1/2) | h\in\h \}$. Note that these $\h(-1/2)\otimes T_{\chi}$ are also irreducible $A(V_L^+)$-modules.
\end{coro}
\begin{proof} It follows immediately from Theorem \ref{atl}.
\end{proof}
Next, we set $$\overline{A(V_L^+)}=A(V_L^+)/A^t_L.$$ For $a\in A(V_L^+)$, we will conveniently denote its image in $\overline{A(V_L^+)}$ by $a$. Clearly, $V_L^{T_{\chi,+}}(0)$ is a $\overline{A(V_L^+)}$-module. Here, $T_{\chi}$ is an irreducible $\hat{L}/K$-module with central character $\chi$ such that $\chi(\iota(\kappa))=-1$.

\begin{thm}\label{qavl+} $\overline{A(V_L^+)}$ is semisimple. In fact, it is isomorphic to $\C[\hat{L}/K]/J$ where $J$ is the ideal of $\C[\hat{L}/K]$ generated by $\kappa K+1$. Consequently, 
\begin{eqnarray*}
\{T_{\chi}\ \ &|& \ \ T_{\chi} \text{ are irreducibles }\hat{L}/K-\text{modules}\\
& &\text{ with central character }\chi\text{ such that }\chi(\iota(\kappa))=-1\ \ \}\end{eqnarray*} is the set of all irreducible $\overline{A(V_L^+)}$-modules. Note that these $T_{\chi}$ are also irreducible $A(V_L^+)$-modules.

\end{thm}
\begin{proof} Recall that there are identity maps from $A(M(1)^+)$ into $A(V_L^+)$ and $A(V_{\Z\alpha}^+)$ into $A(V_L^+)$ for $\alpha\in L-\{0\}$. By following the proof of Lemmas 3.8, 3.9 in \cite{y}, we will obtain that every element in $A(M(1)^+)$ is a constant in $\overline{A(V_L^+)}$. In particular, $[\omega_a]=\frac{1}{16}$ and $[H_a]=\frac{9}{128}$ in $\overline{A(V_L^+)}$ for all $a\in \{1,...,d\}$. Moreover, for $\alpha\in L, u\in M(1)^+$, $[E^{\alpha}]$ commutes $[u]$ in $\overline{A(V_L^+)}$, and $\overline{A(V_L^+)}$ is a direct sum of $M(1)(\theta)^+(0)$. Since $A(V_L^+)(\alpha)=span_{\C}\{[u]*[E^{\alpha}]| u\in M(1)^+\}$ for all $\alpha\in L$, we can conclude that $$\overline{A(V_L^+)}=span_{\C}\{[E^{\alpha}]+A^t_L | \alpha\in L\}.$$ Next, for $\alpha\in L-\{0\}$, we set $B_{\alpha}=2^{\<\alpha,\alpha\>-1}E^{\alpha}$ and $B_0=1$. By following the proof of Lemma 3.11 in \cite{y}, we can show that in $\overline{A(V_L^+)}$, $[B_{\alpha}]*[B_{\beta}]=\epsilon(\alpha,\beta)[B_{\alpha+\beta}]$ for $\alpha,\beta\in L$. Let $\phi:\hat{L}\rightarrow \overline{A(V_L^+)}$ be a linear map defined by $\phi(e_{\alpha})=[B_{\alpha}]$ and $\phi(\kappa)=-1+A^t_L$. Then $\phi$ induces an algebra isomorphism from $\C[\hat{L}/K]/J$ onto $\overline{A(V_L^+)}$. This implies that $\overline{A(V_L^+)}$ is a semisimple algebra.
\end{proof}
\begin{coro}\label{mann} Let $M$ be an $A(V_L^+)$-module such that $A^tm=0$ for all $m\in M$. Then $M$ is an $\overline{A(V_L^+)}$-module. Furthermore, $M$ can be rewritten as a direct sum of irreducible $A(V_L^+)$-modules.\end{coro}
\begin{proof} This follows immediately from Theorem \ref{qavl+}.\end{proof}

\begin{lem}\label{semi1} The Zhu algebra $A(V_L^+)$ is semisimple.
\end{lem}
\begin{proof} Let $M$ be an $A(V_L^+)$-module. We will show that $M$ can be rewritten as a direct sum of irreducible $A(V_L^+)$-modules.
\begin{description}
\item[case 1:] $A^tm=0$ for all $m\in M$. 

Then it follows immediately from Corollary \ref{mann} that $M$ is a direct sum of irreducible $A(V_L^+)$-modules.
\item[case 2:] $A^tm\neq 0$ for some $m\in M$.

First, we will show that $M$ contains a simple $A^t_L$-module. Let $m\in M$ such that $A^tm\neq 0$. Clearly, $A^t_L m$ is an $A^t_L$-module. By Theorem \ref{atl}, we can conclude that $A^t_Lm$ is a direct sum of irreducible $A^t_L$-modules which are also irreducible $A(V_L^+)$-modules. Consequently, $M$ contains a simple $A^t_L$-module. 

Next, we will show that $M$ is semisimple as an $A(V_L^+)$-module. We set $N$ be the direct sum of all irreducible $A^t_L$-submodules of $M$, and we set $$M^0=\{m\in M|A^tm=0\}.$$ Note that $N$ is a sum of irreducible $A(V_L^+)$-modules (cf. Corollary \ref{cirrmod}). Also, $M^0$ is an $A(V_L^+)$-submodule of $M$ since $A^t$ is a 2-sided ideal of $A(M(1)^+)$ and for $m\in M^0$, $v\in A^t$, 
$$[v]\cdot([E^{\alpha}]\cdot m)=([v]*I^t)\cdot([E^{\alpha}]\cdot m)=([v]*[E^{\alpha}])\cdot (I^t\cdot m)=0\text{ for all }\alpha\in L.$$
Let $u\in M-N$. We will show that $u\in M^0$. If $A^tu\neq 0$ then $A^t_L u$ is direct sum of simple $A^t_L$-submodules of $M$. This implies that $A^t_Lu\subset N$ and $u\in N$. This is impossible. Hence $A^tu=0$. Furthermore, $u\in M^0$, and $M-N\subset M^0$. Consequently, $M=N+M^0$. By Corollary \ref{mann}, we can conclude that $M$ is a sum of irreducible $A(V_L^+)$-modules.\end{description}

Consequently, $A(V_L^+)$ is semisimple.

\end{proof}
\begin{coro} Every $A(V_L^+)$-module is semisimple.\end{coro}

\begin{thm} If $L$ is a negative definite even lattice of a finite rank then the vertex algebra $V_L^+$ is rational.
\end{thm}
\begin{proof} It is a consequence of Proposition \ref{pv}, Corollary \ref{pvl+} and  Lemma \ref{semi1}.\end{proof}

\subsection{Case II: $L$ is a non-degenerate even lattice that is neither positive definite nor negative definite}

Following Subsection 3.1, we will prove that the vertex algebra $V_L^+$ is rational when $L$ is a non-degenerate even lattice that is neither positive nor negative definite by showing that the Zhu algebra $A(V_L^+)$ is semisimple. 

For the rest of this Subsection, {\em we assume that $L$ is a non-degenerate even lattice that is neither positive definite nor negative definite}. Next, for every $\alpha\in L$, we set $[B'_{\alpha}]$ in the following way:
\begin{enumerate}
\item if $\alpha\in L$ such that $\<\alpha,\alpha\>\neq 0$, we let 
$$[B'_{\alpha}]=2^{\<\alpha,\alpha\>-1}( I^t*[E^{\alpha}]-\frac{2\<\alpha,\alpha\>}{2\<\alpha,\alpha\>-1}[E^t_{11}]*[E^{\alpha}]).$$ Here, $E^t_{11}$ is defined with respect to an orthonormal basis $\{h_a|1\leq a\leq d\}$ such that $h_1\in\C\alpha$. 
\item if $\alpha\in L-\{0\}$ such that $\<\alpha,\alpha\>=0$, we define \begin{eqnarray*}
[B'_{\alpha}]&=&\frac{1}{2}I^t*[E^{\alpha}]\\
& &+\frac{1}{1-2\<\beta,\beta\>-2\<\gamma,\gamma\>}(\<\gamma,\gamma\>[E^t_{11}]*[E^{\alpha}]+\<\beta,\beta\>[E^t_{22}]*[E^{\alpha}])\\
& &+\frac{\sqrt{\<\beta,\beta\>\<\gamma,\gamma\>}}{1-2\<\beta,\beta\>-2\<\gamma,\gamma\>}([E^t_{12}]*[E^{\alpha}]+[E^t_{21}]*[E^{\alpha}]).
\end{eqnarray*}
Here, $\gamma$, $\beta\in L_{\Q}$ such that $\alpha=\gamma+\beta$, $\<\gamma,\gamma\>>0$, $\<\beta,\beta\><0$, and $\<\gamma,\beta\>=0$. Moreover, $E^t_{11}$, $E^t_{22}$, $E^t_{12}$, $E^t_{21}$ are defined with respect to an orthonormal basis $\{h_a|1\leq a\leq d\}$ so that $h_1\in \C\gamma$, $h_2\in\C\beta$.
\item We set $[B'_0]=[I^t].$
\end{enumerate}
\begin{rem} Let $\alpha\in L$ such that $\<\alpha,\alpha\>=-2k$. It was shown in \cite{j} that $E^{2\alpha}=(1-2k)2^{8k+1}+k2^{8k+6}[\omega_{\alpha}]$. Hence $[B'_{2\alpha}]=1$ on $V_L^{T_{\chi},-}(0)$ for any $\chi$.
\end{rem}

Similar to Subsection 3.1, we have that
for $\alpha\in L$ such that $\<\alpha,\alpha\>\neq 0$, and for $i\in\{1,...,d\}$, 
\begin{enumerate}
\item $[E^t_{ij}]*[B'_{\alpha}]=2^{\<\alpha,\alpha\>-1}[E^t_{ij}]*[E^{\alpha}]$ when $j\in\{2,...,d\}$, and 
\item $[E^t_{i1}]*[B'_{\alpha}]=-\frac{2^{\<\alpha,\alpha\>-1}}{2\<\alpha,\alpha\>-1}[E^t_{i1}]*[E^{\alpha}]$.
\end{enumerate} 
Furthermore, we have the following.
\begin{lem}\label{ebn}
For $1\leq a,b\leq d$, $\alpha\in L$ such that $\<\alpha,\alpha\>\neq 0$, $[E^t_{ab}]*[B'_{\alpha}]=[B'_{\alpha}]*[E^t_{ab}]$.
\end{lem}
Next, we let $$\tilde{A}^t_L=Span_{\C}\{[E^t_{ab}]*[B'_{\alpha}]\ \ |\ \  1\leq a,b\leq d, \ \ \alpha\in L\}.$$ 
Notice that $\tilde{A}^t_L$ is an $A^t$-module. Therefore, $\tilde{A}^t_L$ is a direct sum of $M(1)(\theta)^-(0)$, and for $\alpha\in L$ such that $\<\alpha,\alpha\><0$, we have $[B'_{2\alpha}]=1$ in $\tilde{A}^t_L$.

\begin{lem}\label{ebn2} Let $\alpha\in L-\{0\}$  such that $\<\alpha,\alpha\>=0$. Then
\begin{enumerate}
\item there exists $\beta\in L$ so that $\<\beta,\beta\><0$, $\<\alpha,\beta\><0$ and $[B'_{\alpha}]=[B'_{\alpha+2\beta}]$.
\item For $1\leq a,b\leq d$, we have $[E^t_{ab}]*[B'_{\alpha}]=[B'_{\alpha}]*[E^t_{ab}]$.
\end{enumerate}\end{lem}
\begin{proof} The proof of the first statement is very similar to Lemma 4.7 of \cite{y}. Let $\alpha\in L-\{0\}$  such that $\<\alpha,\alpha\>=0$. By Proposition \ref{info2}, we can conclude that there exist $\beta\in L$, and $u\in A^t$ such that $\<\beta,\beta\><0$, $\<\alpha,\beta\><0$ and $[B'_{\alpha}]=[u]*[B'_{\alpha+2\beta}].$ Since $[B'_{\alpha}]=[B'_{\alpha}]*[B'_{2\beta}]=[B'_{\alpha+2\beta}]$ on $V_L^{T_{\chi,-}}(0)$ for any $\chi$, we can conclude that $[u]=[I^t]$ and $[B'_{\alpha}]=[B'_{\alpha+2\beta}]$ in $\tilde{A}^t_L$. 

2. is a consequence of 1..
\end{proof}

\begin{thm} $\tilde{A}^t_L$ is a 2-sided ideal of $A(V_L^+)$ and $\tilde{A}^t_L\cong A^t\otimes_{\C}\C[\hat{L}/K]/J$ where $\C[\hat{L}/K]$ is the group algebra of $\hat{L}/K$ and $J$ is the ideal of $\C[\hat{L}/K]$ generated by $\kappa K+1$. Furthermore,\begin{eqnarray*}
\{ \ \ \h(-1/2)\otimes T_{\chi}\ \ &|& \ \ T_{\chi} \text{ are irreducibles }\hat{L}/K-\text{modules}\\
& &\text{ with central character }\chi\text{ such that }\chi(\iota(\kappa))=-1\ \ \}\end{eqnarray*} is the set of all irreducible $\tilde{A}^t_L$-modules. Note that these $\h(-1/2)\otimes T_{\chi}$ are also irreducible $A(V_L^+)$-modules.
\end{thm}
\begin{proof} By Proposition \ref{info2}, Lemma \ref{ebn} and Lemma \ref{ebn2}, one can easily show that $\tilde{A}^t_L$ is a 2-sided ideal of $A(V_L^+)$. Next, by following the proof of Theorem \ref{atl} step by step, we will obtain that $\tilde{A}^t_L\cong A^t\otimes_{\C}\C[\hat{L}/K]/J$. The rest of the Theorem is clear.
\end{proof}

Following the Subsection 3.1, we set $$\tilde{A}(V_L^+)=A(V_L^+)/\tilde{A}^t_L.$$ 

\begin{thm} $\tilde{A}(V_L^+)$ is semisimple. In fact, it is isomorphic to $\C[\hat{L}/K]/J$ where $J$ is the ideal of $\C[\hat{L}/K]$ generated by $\kappa K+1$. Consequently, 
\begin{eqnarray*}
\{T_{\chi}\ \ &|& \ \ T_{\chi} \text{ are irreducibles }\hat{L}/K-\text{modules}\\
& &\text{ with central character }\chi\text{ such that }\chi(\iota(\kappa))=-1\ \ \}\end{eqnarray*} is the set of all irreducible $\tilde{A}(V_L^+)$-modules. Note that these $T_{\chi}$ are also irreducible $A(V_L^+)$-modules.
\end{thm}
\begin{proof} The proof is very similar to the proof of Theorem \ref{qavl+}.\end{proof}

\begin{coro}Let $M$ be an $A(V_L^+)$-module such that $A^tm=0$ for all $m\in M$. Then $M$ is an $\tilde{A}(V_L^+)$-module. Furthermore, $M$ can be rewritten as a direct sum of irreducible $A(V_L^+)$-module.\end{coro}

\begin{lem}\label{semi2} The Zhu algebra $A(V_L^+)$ is semisimple.
\end{lem}
\begin{proof} The proof is very similar to Lemma \ref{semi1}.\end{proof}

\begin{thm} If $L$ is a non-degenerate even lattice of a finite rank that is neither positive definite nor negative definite then the vertex algebra $V_L^+$ is rational.
\end{thm}
\begin{proof} It is a consequence of Proposition \ref{pv}, Corollary \ref{pvl+} and  Lemma \ref{semi2}.\end{proof}

\end{document}